\documentclass[12pt,english]{elsarticle}
\usepackage[T1]{fontenc}
\usepackage[latin9]{inputenc}
\usepackage{amsmath}
\usepackage{amssymb}
\usepackage{graphicx}
\usepackage{esint}

\makeatletter
\usepackage{url}

\makeatother

\usepackage{babel}
\begin{document}

\title{On the error propagation of semi-Lagrange and Fourier methods for
advection problems\tnoteref{label1}}
	\tnotetext[label1]{This work is supported by the Austrian Science Fund (FWF) -- project id: P25346.}

	\author[uibk]{Lukas Einkemmer\corref{cor1}}
	\ead{lukas.einkemmer@uibk.ac.at}

	\author[uibk]{Alexander Ostermann}
	\ead{alexander.ostermann@uibk.ac.at}

	\address[uibk]{Department of Mathematics, University of Innsbruck, Austria}
	\cortext[cor1]{Corresponding author}

\begin{abstract}
In this paper we study the error propagation of numerical schemes
for the advection equation in the case where high precision is desired.
The numerical methods considered are based on the fast Fourier transform,
polynomial interpolation (semi-Lagrangian methods using a Lagrange
or spline interpolation), and a discontinuous Galerkin semi-Lagrangian
approach (which is conservative and has to store more than a single
value per cell). 

We demonstrate, by carrying out numerical experiments, that the worst
case error estimates given in the literature provide a good explanation
for the error propagation of the interpolation-based semi-Lagrangian
methods. For the discontinuous Galerkin semi-Lagrangian method, however,
we find that the characteristic property of semi-Lagrangian error
estimates (namely the fact that the error increases proportionally to
the number of time steps) is not observed. We provide an explanation
for this behavior and conduct numerical simulations that corroborate
the different qualitative features of the error in the two respective
types of semi-Lagrangian methods.

The method based on the fast Fourier transform is exact but, due to
round-off errors, susceptible to a linear increase of the error in
the number of time steps. We show how to modify the Cooley--Tukey
algorithm in order to obtain an error growth that is proportional
to the square root of the number of time steps. 

Finally, we show, for a simple model, that our conclusions hold true
if the advection solver is used as part of a splitting scheme.
\end{abstract}
\maketitle

\section{Introduction}

The accurate numerical simulation of advection dominated
problems is an important problem in many scientific applications.
However, due to the non-dissipative nature of the equations considered,
care has to be taken to obtain a stable numerical scheme (see, for
example, \citep{morton1988}). A large body of research has been accumulated
that describes finite difference, finite volume, and finite element
discretizations of such problems. However, for advection-dominated
problems so called semi-Lagrangian methods offer a competitive alternative.
These methods integrate the characteristics back in time and consequently
have to use some interpolation scheme to reconstruct the desired value
at the grid points. Strictly speaking, semi-Lagrangian methods can
only be applied to systems of first-order differential equations.
However, in many instances, first-order systems arise from the splitting
of more complicated equations or constitute the linear part of an
evolution equation (which is then treated separately from the nonlinearity).
Consequently, semi-Lagrangian methods have been used extensively in
applications ranging from fluid dynamics to plasma physics (see e.g.
\citep{staniforth1991} and \citep{filbet2003}). Such an approach
is especially promising, if the characteristics (of a sub-problem)
can be computed analytically; this can be done, for example, in context
of the Vlasov--Poisson equations. In addition, semi-Lagrangian methods
do not impose a Courant--Friedrichs--Lewy (CFL) condition.

In some problems, methods based on the FFT (fast Fourier transform)
can also be employed; this is the case for tensor product domains
(see e.g. \citep{klimas1994}). Compared to FFT based methods, the
semi-Lagrangian methods provide a local approximation (which is important
in the context of parallelization). They are more easily applicable
to non-periodic boundary conditions (due to the absence of Gibbs'
phenomenon) and usually are better suited to handle nonlinearities.
Using the fast Fourier transform, on the other hand, allows us to
solve the linear advection equation exactly (in infinite precision
arithmetics).

In most scientific applications a tolerance of say $10^{-3}$ is sufficient.
In this case, the main research goal is to construct more efficient
algorithms and to implement better step size control mechanism. Also
methods that preserve certain invariants of the continuous system
are of interest in that context. 

However, a number of applications have been identified where double
precision floating point numbers are not sufficient. A proposed remedy
is to (selectively) use 128-bit floating point numbers. Note, however,
that this procedure is accompanied by a significant reduction in performance.
Examples of such problems range from the investigation of vortex sheet
roll-ups in fluid dynamics to electromagnetic scattering phenomena
(for an excellent review article see \citep{bailey2005}). 

Also, concern has been raised in recent years with regard to the reproducibility
of numerical simulations; especially if such simulations are conducted
on different computer systems. In \citep{he2000}, for example, it
is demonstrated that climate codes show significantly different result
depending on the number of processors that are employed in the simulation.
One popular choice of numerical methods for atmospheric modeling are
semi-Lagrangian methods. 

Furthermore, due to the diminishing gain in per core CPU (central
processing unit) performance, massively parallel computing architectures,
such as GPUs and the Xeon Phi, have become more and more common. Usually
in such situations the memory per core is significantly smaller than
in more traditional cluster systems. In addition, single precision
floating point performance is usually faster than double precision
floating point performance (for example, the CUDA FFT single precision
implementation achieves a speedup of about $2.5$ compared to the
double precision implementation \citep{cufft}). Such considerations
make the use of single precision floating point numbers attractive
for some applications.

Therefore, our goal in this paper is to study the error propagation
in a context where results close to machine precision are of interest
or where the error has to be tightly controlled. We will limit ourselves
to the advection equation (on a finite spatial interval)
\begin{equation}
\partial_{t}u(t,x)+v\partial_{x}u(t,x)=0.\label{eq:advection}
\end{equation}

For this model problem, we consider the time evolution of the interplay
of round-off and discretization errors for semi-Lagrangian and FFT
based methods. The discretization of (\ref{eq:advection}) is important
as it is a building block for many more involved schemes (for example,
in the context of splitting methods or in methods where a linear part
is treated differently). Such schemes are applied, for example, in
fluid dynamics or to solve the Vlasov equation. 

For future reference, we note that the analytic solution of (\ref{eq:advection})
can be easily written down as
\[
u(t,x)=u(0,x-vt),
\]
where in this paper we always assume that $v$ is a constant independent
of $x$ and $t$.

\section{Description and error bounds}

In this section we will describe how to use the semi-Lagrangian method
(section \ref{sub:semi-Lagrange}) and the fast Fourier transform
(section \ref{sub:FFT}). Furthermore, we will discuss some theoretical
results concerning the discretization error of these schemes.

\subsection{The semi-Lagrangian method\label{sub:semi-Lagrange}}

A time step in the semi-Lagrangian method for the $i$th grid point
is computed as follows:
\[
u_{n}\left(x_{i}\right)=u_{n-1}\left(X_{x_{i}}(\tau)\right),
\]
where $u_{n}$ is the numerical solution after $n$ time steps, $\tau$
is the time step size, and the characteristics of equation (\ref{eq:advection})
are given by $X_{x}(\tau)=x-v\tau$. Note that $X_{x_{i}}(\tau)$
does not necessarily coincide with any grid point. Therefore, a sufficiently
accurate interpolation procedure has to be used in order to extend
the values stored at the grid to the entire domain. Both (continuous)
spline interpolation as well as (discontinuous) Legendre or Lagrange
interpolation are popular choices. Furthermore, similar to discontinuous
Galerkin methods, multiple coefficients can be stored for each cell
(which yields a local reconstruction at the expense of additional
memory demands).

It is well known (see e.g. \citep{besse2004} and \citep{einkemmer2014})
that in the case of semi-Lagrangian methods the error of the fully
discretized problem can be estimated by
\begin{equation}
\Vert u(n\tau,x)-u_{n}(x)\Vert\leq C\left(h^{q}+\frac{h^{q}}{\tau}\right),\label{eq:error-estimate-sL}
\end{equation}
where $h$ is the size of the spatial grid, and $q$ is the order
of the space discretization (for example, the order of the Lagrange
interpolation). Such estimates are usually derived in one (or all)
of the $L^{p}$ norms. Note that there is no term proportional to
$\tau$ as the characteristics are known analytically. However, results
can be derived that take time discretization errors into account (see
e.g. \citep{besse2004} and \citep{einkemmer2014}).

The first term on the right-hand side in equation (\ref{eq:error-estimate-sL})
is the usual spatial error term. The disturbing implication, however,
is that the second term is directly proportional to the number of
time steps taken. This, however, is to be expected, as in each projection
an error proportional to $h^{q}$ is made. In the worst case these
errors accumulate to give the above mentioned bound. Note that the
second term is usually not present if Eulerian discretization methods
(for example finite-difference schemes) are employed.

Similarly how the discontinuous Galerkin (dG) method considers a
polynomial reconstruction that only uses data attached to the corresponding
cell, it is possible to formulate a semi-Lagrangian method based on
this reconstruction. This dG semi-Lagrangian method also requires
only the data from neighboring cells in order to perform a time step
(for arbitrary order in space). The scheme is described in some detail
in \citep{einkemmer2014-2QOYZWMPACZAJ2MABGMOZ6CCPY} and \citep{einkemmer2014}.
Since there is no ambiguity, in this paper, we will refer to this
scheme as the dG approximation.

Let us duly note that the second term in equation (\ref{eq:error-estimate-sL})
is a worst case estimate that is valid for all the interpolation and
projection schemes discussed here. Thus, in the next section we will
investigate the actual numerical behavior for the spline interpolation,
the Lagrange interpolation, and the discontinuous Galerkin approximation.

\subsection{The fast Fourier transform\label{sub:FFT}}

In case of the FFT based scheme we compute the Fourier transform of
(\ref{eq:advection}) which gives
\begin{equation}
\partial_{t}\hat{u}(t,k)+ikv\hat{u}(t,k)=0,\label{eq:fourier-advec-eq}
\end{equation}
where we have denoted the Fourier component with frequency $k\in\mathbb{Z}$
by $\hat{u}(t,k)$. We can easily solve (\ref{eq:fourier-advec-eq})
to get
\[
\hat{u}(t,k)=\mathrm{e}^{-ikvt}\hat{u}(0,k).
\]
Therefore, the advection in Fourier space is described by the multiplication
with an appropriate phase factor. The numerical scheme truncates the
Fourier series. Therefore, we only consider $-m\leq k\leq m$. Since
$u(t,k)$ is a real function, only the non-negative frequencies have
to be stored in memory. 

If exact arithmetics is used, we obtain an error bound that depends
on the spatial regularity of the solution. For $u(t,\cdot)\in\mathcal{C}^{m}$
we have (as a consequence of the Riemann--Lebesgue lemma)
\[
\Vert u(n\tau,x)-u_{n}(x)\Vert\leq Ch^{m-1},
\]
where $h$ is the grid size. In particular, for $u(t,\cdot)\in\mathcal{C}^{\infty}$
the convergence is super-polynomial in the number of grid points.
This is the only discretization error provided that exact arithmetics
is employed.

\section{Numerical investigation\label{sec:Numerical-investigation}}

The purpose of this section is to present the results from a number
of numerical simulations conducted. It will soon be apparent that
the actual computations displays a more complicated behavior as would
be expected from the error estimates discussed in the previous section.

\begin{figure}
\noindent \begin{centering}
\includegraphics[width=12cm]{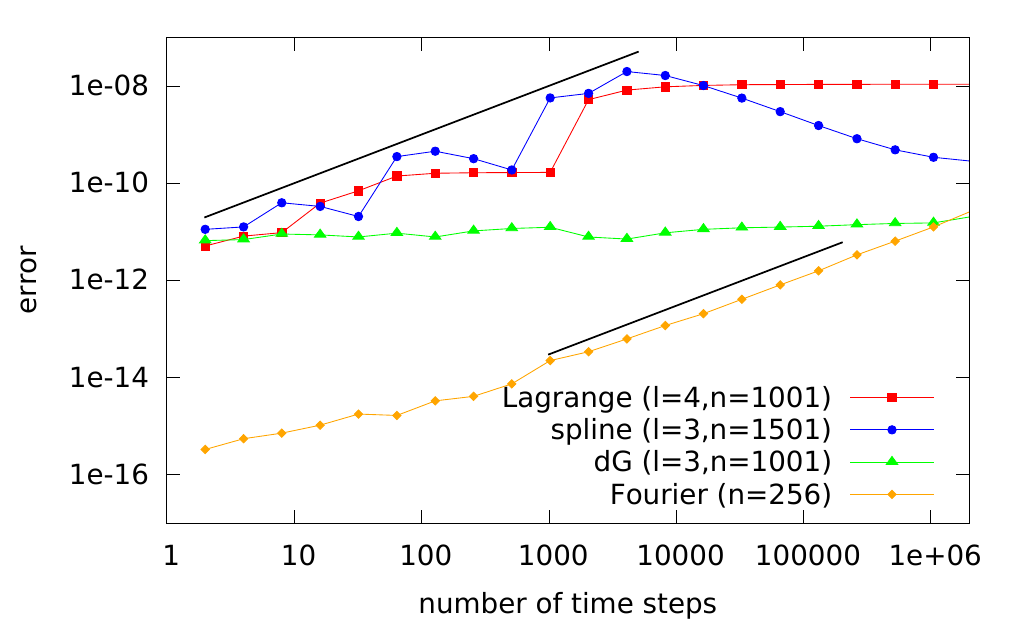}
\par\end{centering}

\caption{The $L^{\infty}$ error of the Lagrange, spline, dG, and FFT based
methods, as a function of the number of time steps, is shown. The
FFT routine from the FFTW library is used. The polynomial degree is
denoted by $l$ and the number of cells/grid points is denoted by
$n$. As a reference two black lines of slope $1$ are drawn. \label{fig:comparison-all}}
\end{figure}

In Figure \ref{fig:comparison-all}, we compare the error propagation
for a Lagrange interpolation, a discontinuous Galerkin method, and
the Fourier approximation using the initial value
\begin{equation}
u(0,x)=\frac{1}{2+\cos\pi x}\label{eq:1overcos-iv}
\end{equation}
on the interval $[-1,1]$ with periodic boundary conditions (the same
interval and periodic boundary conditions are used for all simulations
in this paper). As expected, initially the FFT method achieves a performance
close to machine precision. Note, however, that the error growth is
linear in the number of time steps. However, from a stochastic description
of the round-off error one would expect an error growth proportional
to the square root in the number of time steps.  Let us postpone
the detailed investigation of this issue until section \ref{sec:fft-round-off-error}.

Furthermore, even though the worst case error estimate for a general
semi-Lagrangian method does include the term proportional to the number
of steps, this is only observed for the Lagrange and the spline interpolation
(see Figure \ref{fig:comparison-all}). However, the dG method does
not exhibit such a behavior. In fact, there is almost no error propagation
even after more than $10^{6}$ steps in time have been conduced. In
this instance, an error in the initial value that is orders of magnitude
away from machine precision can still be competitive with the FFT
based scheme (which shows a linear propagation of the error).

For both the Lagrange and the spline interpolation some oscillations
do occur. In general, however, a linear error propagation is observed.

\begin{figure}
\noindent \begin{centering}
\includegraphics[width=12cm]{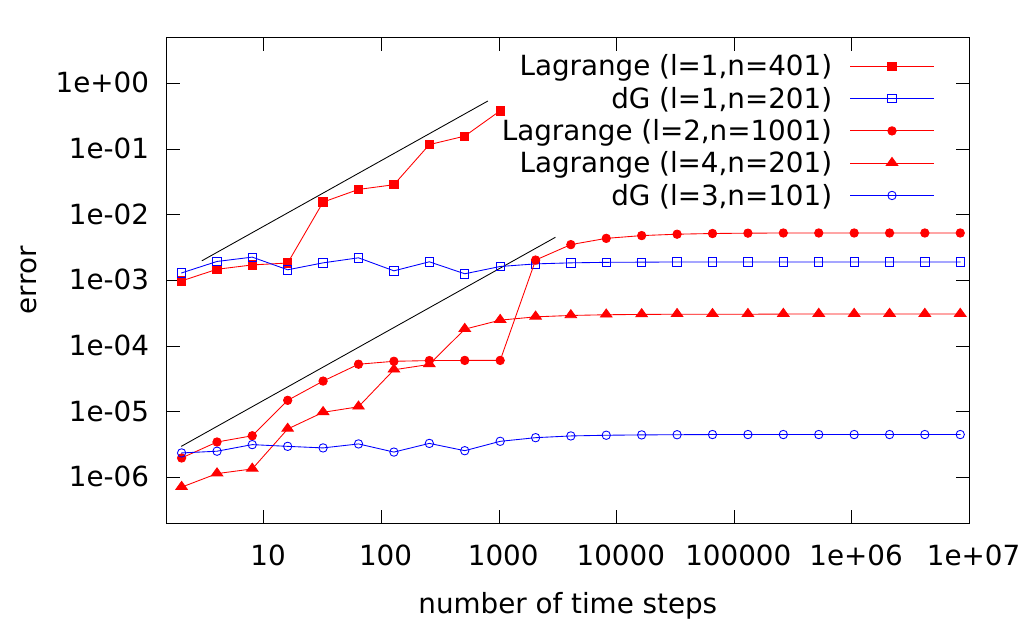}
\par\end{centering}

\caption{The $L^{\infty}$ error as a function of the number of time steps
is shown for a number of different configurations (the interpolation
degree $l$ and the number of grid points/cells $n$ of the Lagrange
and dG interpolation methods are varied). The initial value given
in (\ref{eq:cosine-iv}) is used. As a reference two black lines of
slope $1$ are shown. \label{fig:different-order}}
\end{figure}

In Figure \ref{fig:different-order} we compare different polynomial
degrees for both the Lagrange and dG based methods with the initial
value
\begin{equation}
u(0,x)=\cos4\pi x.\label{eq:cosine-iv}
\end{equation}
As before, we initially observe a linear error growth, which yields
a reduction in accuracy of at least three orders of magnitude for
the Lagrange interpolation. No such behavior is observed for the dG
method.

\section{A theoretical investigation of the semi-Lagrangian methods\label{sec:theoretical-investigation}}

In the previous section we observed a remarkable difference in the
error propagation for the semi-Lagrangian methods employing Lagrange
or spline interpolation on the one hand and the discontinuous Galerkin
method on the other hand. The conjecture is that in case of the discontinuous
Galerkin scheme the errors made in each step average to zero over
the course of many time steps, while for the interpolation methods
they are in fact well described by the worst case estimate stated
in the introduction.

It is our goal now to provide an explanation for this behavior in
the case of the linear Lagrange interpolation. For simplicity, let
us drop the time dependence of $u$ in this section. Thus, we denote
the initial value for a given time step by $u(x)$ and the $i$th
grid point by $x_{i}$. A single time step of size $\tau$ (and with
$v=1$) then gives the new value at the $i$th grid point (which we
denote by $\tilde{u}(x_{i})$)
\begin{equation}
\tilde{u}(x_{i})=\frac{1}{h}\left(\tau u(x_{i-1})+(h-\tau)u(x_{i})\right).\label{eq:utilde}
\end{equation}
Note that we have restricted ourselves to $\tau<h$. This is justified
as each advection with $\tau>h$ can be decomposed into an advection
that can be solved exactly, where $\tau$ is a multiple of $h$, and
an advection for which $\tau<h$ holds true.

The corresponding extension to the entire domain is then given by
(where $\xi\in[0,h]$ is restricted to the cell under consideration)
\begin{equation}
\tilde{u}(x_{i}+\xi)=\frac{1}{h}\left((h-\xi)\tilde{u}(x_{i})+\xi\tilde{u}(x_{i+1})\right).\label{eq:utilde-everywhere}
\end{equation}

Before we proceed, let us make two remarks. First, in what follows,
we compute the error as compared to the advected initial value 
which lies in the space of piecewise linear functions (and not to the analytic
initial value). This is in fact the correct choice as we are interested
in the error propagation and not in the initial projection error (which
clearly is bounded by $Ch^{2}$ in this case). Second, the error is
a function of two variables; the position $x=x_{i}+\xi$ and the size
of the time step $\tau$.

Let us consider the average (with respect to the spatial variable
in a single cell, which is denoted by $\xi$) for the exact advection
\begin{align*}
\frac{1}{h}\int_{0}^{h}u(x_{i}+\xi-\tau)\,\mathrm{d}\xi & =\frac{1}{h}\int_{h-\tau}^{h}u(x_{i-1}+\xi)\,\mathrm{d}\xi+\frac{1}{h}\int_{0}^{h-\tau}u(x_{i}+\xi)\,\mathrm{d}\xi\\
 & =\frac{\tau^{2}}{2h^{2}}u(x_{i-1})+\frac{(h+\tau)^{2}-3\tau^{2}}{2h^{2}}u(x_{i})+\frac{(h-\tau)^{2}}{2h^{2}}u(x_{i+1})
\end{align*} and for the linear Lagrangian interpolation (using (\ref{eq:utilde-everywhere})
and (\ref{eq:utilde}))
\begin{align*}
\frac{1}{h}\int_{0}^{h}\tilde{u}(x_{i}+\xi)\,\mathrm{d}\xi & =\frac{\tau}{2h}u(x_{i-1})+\frac{1}{2}u(x_{i})+\frac{h-\tau}{2h}u(x_{i+1}).
\end{align*}

Now, since we are interested in simulations where a large number of
steps has to be conducted it is reasonable to assume that for a fixed
cell in the computational grid the step size $\tau\in[0,h]$ is uniformly
distributed across that interval. Thus, the average in a single cell
(averaged in both space as well as step size) for the exact solution
is given by
\[
\frac{1}{h}\int_{0}^{h}\frac{1}{h}\int_{0}^{h}u(x_{i}+\xi-\tau)\,\mathrm{d}\xi\mathrm{d}\tau=\frac{1}{6}\left(u(x_{i-1})+4u(x_{i})+u(x_{i+1})\right)
\]
and for the linear Lagrangian interpolation we get
\[
\frac{1}{h}\int_{0}^{h}\frac{1}{h}\int_{0}^{h}\tilde{u}(x_{i}+\xi)\,\mathrm{d}\xi\mathrm{d}\tau=\frac{1}{4}\left(u(x_{i-1})+2u(x_{i})+u(x_{i+1})\right).
\]
The double averaged error $\overline{\overline{e}}$ in a single cell
is therefore given by
\[
\overline{\overline{e}}=\frac{1}{12}\left(u(x_{i-1})-2u(x_{i})+u(x_{i+1})\right)\approx\frac{h^{2}}{12}u^{\prime\prime}(x_{i}).
\]

Since there is a non-zero average error this error is amplified in
each time step and gives a linear error propagation. Similarly to
this computation we can also derive results for higher degree Lagrange
interpolation. For the case of piecewise quadratic polynomials, for
example, we get
\begin{align*}
\overline{\overline{e}} & =\frac{1}{144}\left(u(x_{i-2})+2u(x_{i-1})-12u(x_{i})+14u(x_{i+1})-5u(x_{i+2})\right)\\
 & \approx-\frac{h^{3}}{24}u^{(3)}(x_{i}).
\end{align*}

In contrast, the discontinuous Galerkin method by construction preserves
the average exactly. This is a necessary condition (but not a sufficient
one) in order for the errors to cancel out on average. This is, of
course, a different way of stating that the dG scheme is locally conservative
for any step size while the Lagrange interpolation is not. The fact
that most semi-Lagrangian schemes are not conservative is well established
in the literature. To remedy this deficiency usually high-order methods
are employed which provide sufficient accuracy to keep the violation
in mass conservation to an acceptable level. 

However, in our context this is not a remedy since, as we have observed
in the previous section, even a Lagrange interpolation that is of
high accuracy will loose at least three orders of magnitude in precision
after approximately $10^{3}$ time steps. This is no concern if an
approximation correct to three digits is desired. However, it is a
significant drawback if we are interested in accuracies close to machine
precision.

Note that since the average error is proportional to the second derivative
we would expect that this behavior can be observed in numerical simulations.
To that end we consider the convex initial value
\begin{equation}
u(0,x)=(x-1)(x+1)\label{eq:convex-iv}
\end{equation}
and the concave initial value
\begin{equation}
u(0,x)=-(x-1)(x+1).\label{eq:concave}
\end{equation}

\begin{figure}
\noindent \begin{centering}
\includegraphics[width=12cm]{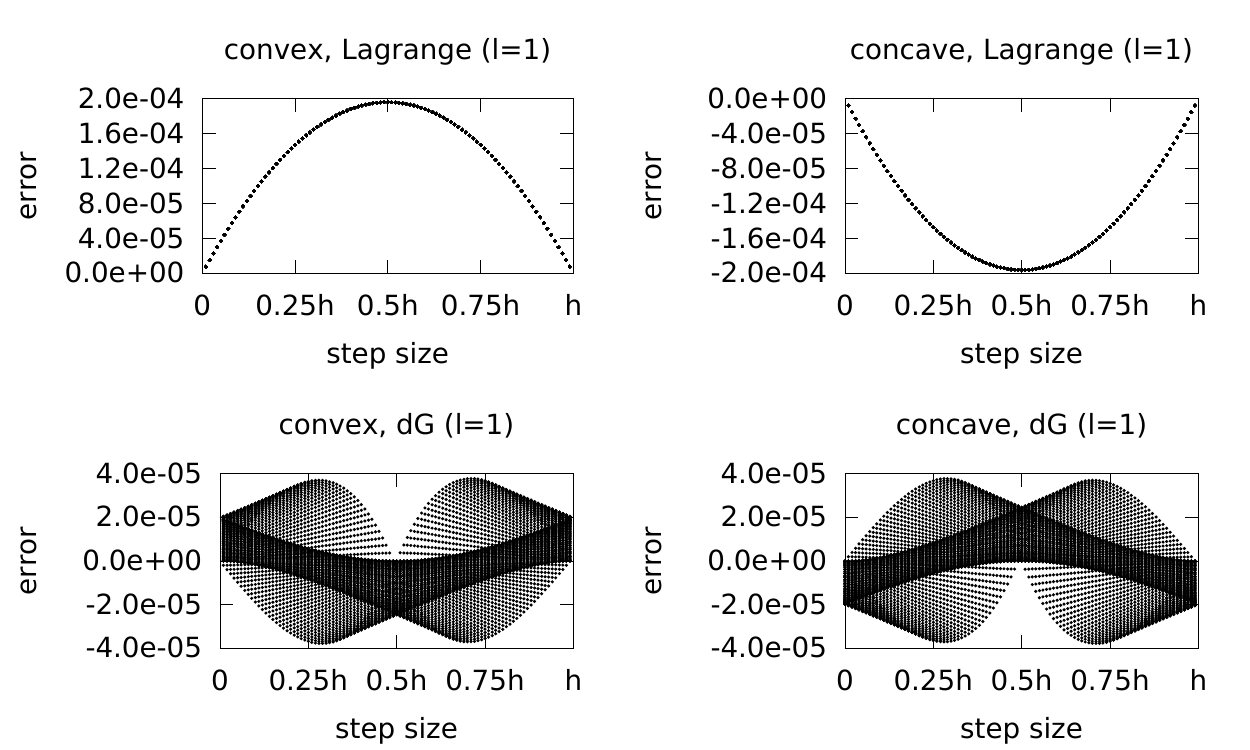}
\par\end{centering}

\caption{The error as a function of the step size $\tau\in[0,h]$ is shown
for the convex (\ref{eq:convex-iv}) and concave (\ref{eq:concave})
initial values. In case of the Lagrangian interpolation the error
at the different grid points is almost identical and thus the points
fall on a line. For the dG method the error computed on different
positions (within a single cell) are displayed at the corresponding
$\tau$ value. \label{fig:convex-concave}}
\end{figure}

In fact Figure \ref{fig:convex-concave} shows the expected behavior.
Note that the errors made do depend on both the step size $\tau$
as well as the position $x$. For the Lagrange interpolation we have
plotted the error at the different grid points at the same $\tau$
value (these points are almost identical and thus indistinguishable
in the plot). In the case of the discontinuous Galerkin method we
want to demonstrate that the error inside a single cell cancels out.
Therefore, we have plotted the error at different points (but inside
the same cell) at the same $\tau$ value.

\begin{figure}
\noindent \begin{centering}
\includegraphics{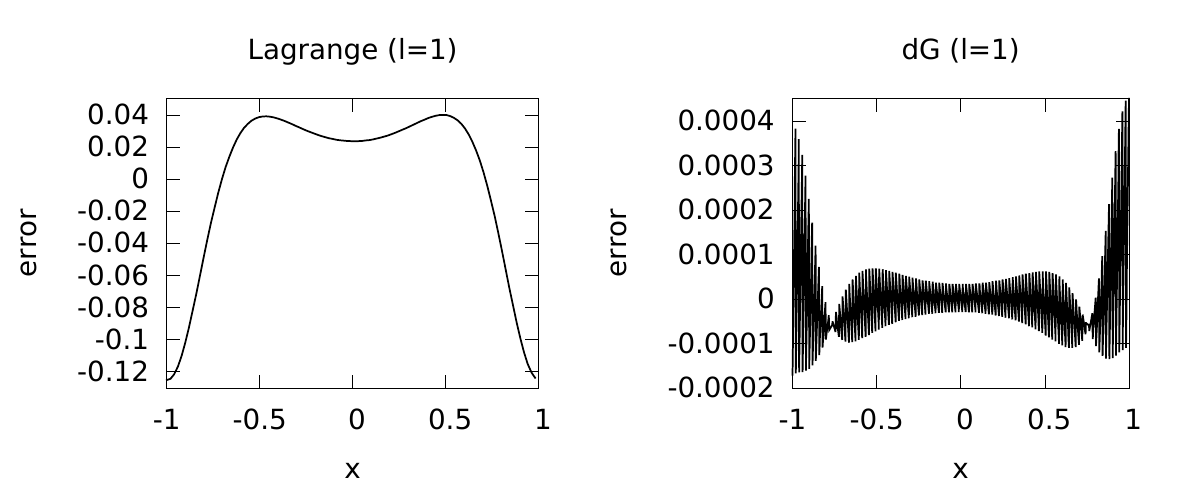}
\par\end{centering}

\caption{The pointwise error after $10^{4}$ steps of size $\tau=2\cdot10^{-4}$
for the Lagrange and dG method (with $101$ grid points) is shown.
The left plot (Lagrange interpolation) matches the second derivative
of (\ref{eq:1overcos-iv}) (the initial value used in this simulation).
In case of the discontinuous Galerkin approximation the error displays
small oscillations of high frequency. \label{fig:pointwise-error}}
\end{figure}

In addition, we show the error as a function of the spatial variable
in Figure~\ref{fig:pointwise-error}. As would be expected we observe
an error that is similar to the second derivative of (\ref{eq:1overcos-iv})
for the Lagrange interpolation and an oscillatory error for the dG
method.

\section{Fast Fourier transform round-off error propagation\label{sec:fft-round-off-error}}

\begin{figure}
\noindent \begin{centering}
\includegraphics[width=12cm]{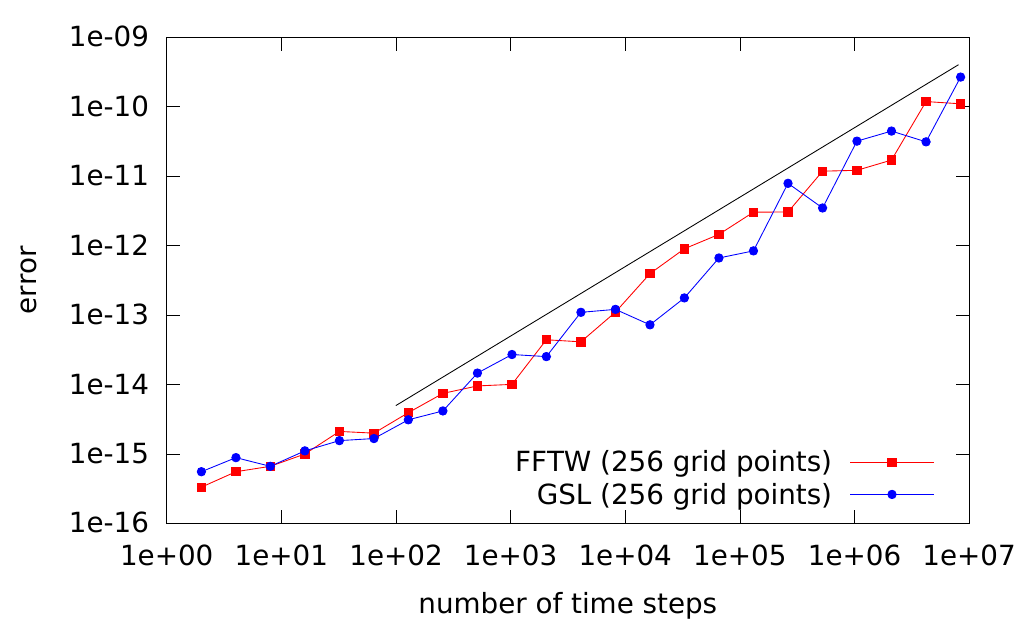}
\par\end{centering}

\caption{The $L^{\infty}$ error of the FFT based advection as a function of
the number of time steps. Results for the FFTW and GSL library are
shown. As a reference a black line with slope $1$ is also displayed.
\label{fig:fft-double}}
\end{figure}

In the numerical simulations conducted in section \ref{sec:Numerical-investigation}
we observed that the error growth is linear in the number of time
steps (and not proportional to the square root as we would expect
from a pure propagation of round-off errors). This behavior is consistent
across numerical libraries; in Figure~\ref{fig:fft-double} results
using the FFTW (Fastest Fourier Transform in the West) library and
the radix-2 implementation found in GSL (GNU Scientific Library) are
shown. 

\begin{figure}
\noindent \begin{centering}
\includegraphics[width=12cm]{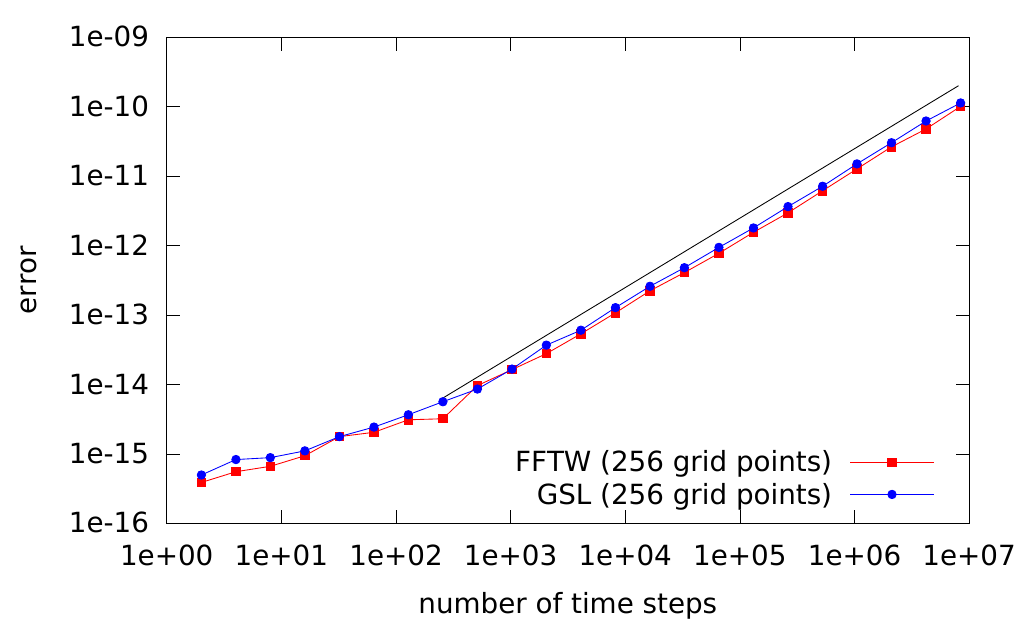}
\par\end{centering}

\caption{$L^{\infty}$ error of the FFT based advection as a function of the
number of time steps. The phase factor needed for the translation
and the corresponding combination with the signal is computed in high
precision arithmetics and then rounded down to double precision. As
a reference a black line with slope $1$ is also shown. \label{fig:fft-phase-factor}}
\end{figure}

An obvious explanation is to suggest that the phase factor used to
compute the advection incurs some additional round-off errors. This
phenomenon is well known in the context of reducing the round-off
errors introduced by FFT routines (see e.g.~\citep{schatzman1996}).
In that context, care has to be taken that the twiddle factors are
computed to sufficient accuracy. A similar approach can be used to
compute the phase factor in the advection. However, Figure \ref{fig:fft-phase-factor}
clearly shows that the error growth is still linear in the number
of time steps.

Note that the FFT algorithm as originally proposed by Cooley and Tukey
is usually not implemented in high-performance FFT libraries (such
as FFTW). A number of additional optimizations are performed. For
the FFTW library a discussion can be found in \citep{fftwaccuracy}.
These optimizations often have a significant impact on accuracy. Therefore,
accuracy benchmarks are performed in order to verify that the round-off
errors are reasonable (see, for example, \citep{fftwaccuracy}). However,
in this paper we are interested in error propagation for a large number
of time steps and not primarily with the scaling of the round-off error
for large problem sizes. 

\begin{figure}
\noindent \begin{centering}
\includegraphics[width=12cm]{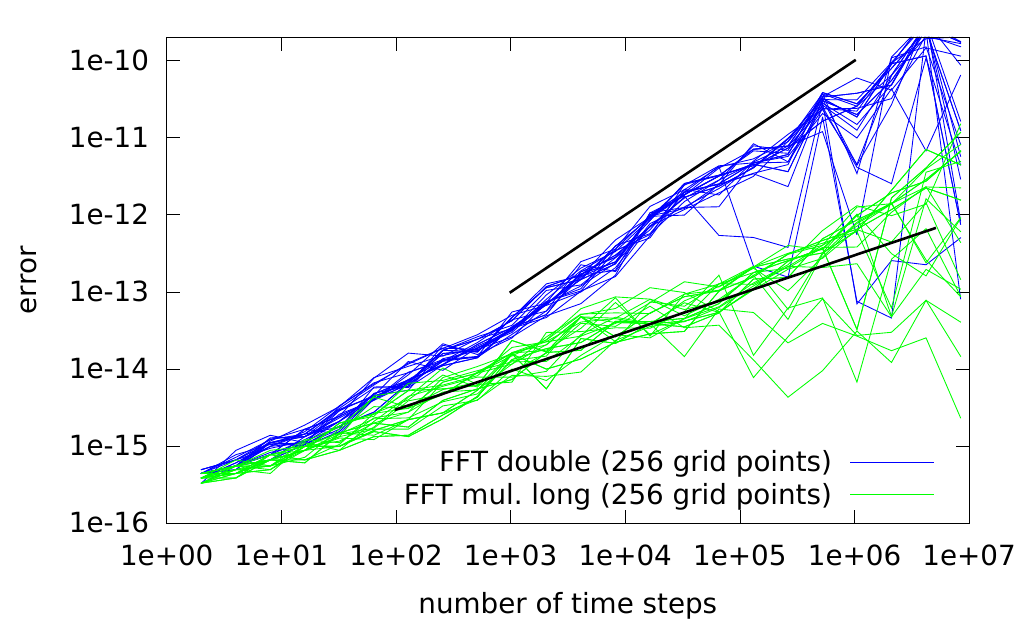}
\par\end{centering}

\caption{The $L^{\infty}$ error of the FFT based advection as a function of
the number of time steps is shown for a number of different initial
values as given in equation (\ref{eq:random-initial}). The plain
FFT implementation following the Cooley--Tukey algorithm (blue) and
a plain FFT algorithm where the multiplications for the phase factor
as well as in the Cooley--Tukey algorithm are carried out in 80-bit
arithmetics (green) are compared. Note that in the latter case the
storage requirement is still essentially the same (that is we do only
have to store the input and output vector, in double precision, as
well as the twiddle factors). As a reference two black lines with
slope $1$ (upper line) and slope $1/2$ (lower line) are shown. \label{fig:fft-80bit}}
\end{figure}

In Figure \ref{fig:fft-80bit} the error propagation for the plain
Cooley--Tukey algorithm is shown for the initial value
\begin{equation}
u(0,x)=\frac{1}{2+\cos(\pi x+\varphi)},\label{eq:random-initial}
\end{equation}
where $\varphi$ is chosen at random. Even though the error propagation
is significantly reduced compared to the FFTW implementation, at least
four orders of magnitude are lost and the initial error growth is
still linear.

If, in addition, the multiplications in the fast Fourier transform
are computed to a higher precision, the error growth shows a behavior
that is roughly proportional to the square root of the number of steps.
In Figure \ref{fig:fft-80bit} the multiplications have been implemented
in the 80bit extended precision type%
\footnote{This data type is usually mapped to \texttt{long double} by the C/C++
compiler. Note, however, that this is not required by the C/C++ standard
and is not entirely consistent across different compiler implementations.
In all our studies we have used the GNU Compiler Collection (versions
4.6 and 4.7).%
} implemented in the x86 hardware.

In the naive implementation used here we observe a reduction in performance
of about 10-15\%. Note, however, that if vectorization is used such
procedures can result in a more severe reduction in performance. 

On the other hand, this approach would also be advantageous, if double
precision floating point computations significantly impact performance
(such as commonly found on GPU systems). In this case the input and
the output would be stored as single precision floating point numbers
and multiplications would be performed as double precision (with appropriately
computed twiddle factors) and then correctly rounded to single precision.

\section{Splitting of an advection equation with source term}

\begin{figure}
\noindent \begin{centering}
\includegraphics[width=12cm]{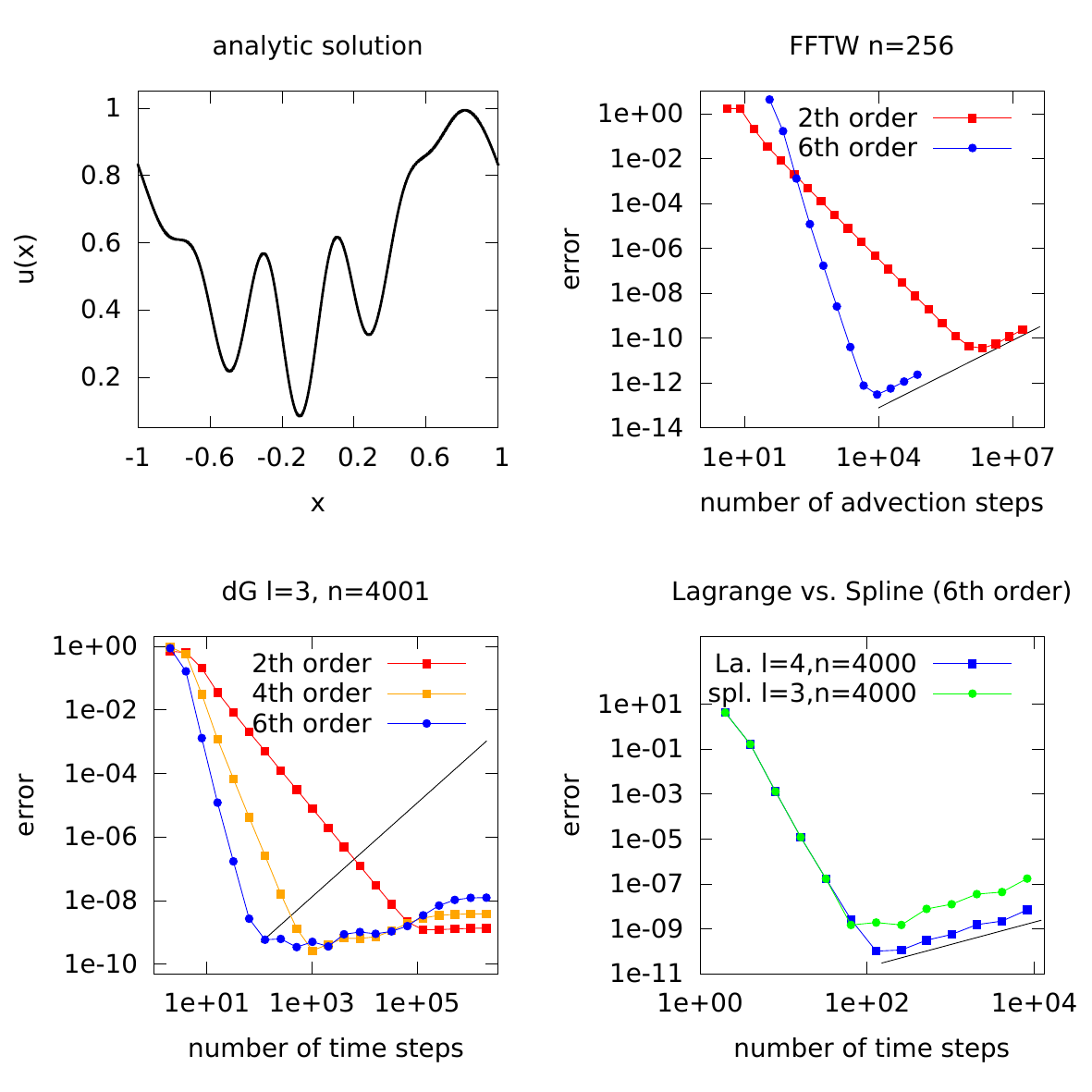}
\par\end{centering}

\caption{The solution of (\ref{eq:advection-source}) is computed for the final
time $t=1.8$. The analytical solution is shown on the top left. The
numerical solution is computed using a splitting approach. The results
for the FFT based method (top right), the discontinuous Galerkin (dG)
method (bottom left), and the Lagrange and spline interpolations (bottom
right) are shown. The black lines have slope $1$ and are shown for comparison. \label{fig:advection-source}}
\end{figure}

In this section we consider the advection equation supplemented by
a (position dependent) source term. That is, we consider
\begin{equation}
\partial_{t}u(t,x)+v\partial_{x}u(t,x)=s(x),\label{eq:advection-source}
\end{equation}
where the source term is chosen as $s(x)=(1+\cos\pi x)\cos5\pi x$
and periodic boundary conditions are imposed. The solution (for the
initial value given in (\ref{eq:1overcos-iv})) can be easily computed
and is plotted in Figure~\ref{fig:advection-source} (top left).

For the time integration we employ the second order Strang splitting
scheme as well as the 6th order scheme constructed from it by composition
(see, for example, \citep{hairer2006}). An approximation to the solution
of the first sub-problem (the advection equation) can be computed
as described in the previous sections. The remaining sub-problem (corresponding
to the source term) is easily solved analytically. Note, however,
that in the case of the dG method only the coefficients (in the Legendre
expansion) are stored in memory. Thus, in each cell, we evaluate the
approximation at the Gauss--Legendre points. We then compute the solution
of the sub-problem on these points and use the result to reconstruct
the coefficients.

The numerical simulations conducted confirm the observations made
in the previous sections. For the FFT based method (using the FFTW
library) we see, after a decrease in the error due to the time discretization
error, the characteristic linear increase in the error (see Figure
\ref{fig:advection-source} top right). Furthermore, we observe that,
as expected, the lowest error we can achieve for a given time discretization
is only dependent on the number of advection steps that have to be
made in order to reach that tolerance (thus the method of order six
is clearly the preferred choice in this case). All these phenomena
are due to round-off errors only.

The spline and Lagrange interpolations show a similar behavior. However,
in this case the linear increase in the error (for further decreasing
step size) is dependent on the space discretization (see Figure \ref{fig:advection-source}
bottom right). 

For the dG method, on the other hand, the minimal error that can be
achieved (for a given space discretization) is almost independent
of the numerical method used in time. Of course, the 6th order method
is usually more efficient if high precision is desired. However, if
a large number of time steps are taken with the Strang splitting scheme,
a similar accuracy than for the 6th order method can be achieved (without
changing the space discretization). Furthermore, there is no linear
increase in the error. All these observations are in line with the
observations made in sections \ref{sec:Numerical-investigation} and~\ref{sec:theoretical-investigation}.

\section{Conclusion}

In this paper we have considered the error propagation for the advection
equation in the case where high precision is desired. The numerical
methods considered exhibit a variety of different phenomena. 

In case of the fast Fourier transform method round-off errors are
the primary concern. A number of libraries that implement the FFT
show a linear error growth in the number of time steps. However, if
the multiplication of the Fourier coefficients with the twiddle factors
is performed to sufficient accuracy the growth in the error is only
proportional to the square root in the number of time steps. 

Furthermore, we have shown that the term proportional to the number
of time steps, that is routinely obtained in error estimates for semi-Lagrangian
methods, is not observed for all semi-Lagrangian schemes. In fact
it is true that the qualitative features of the error is markedly
different for the interpolation (Lagrange as well as spline) and the
discontinuous Galerkin based semi-Lagrangian schemes considered in
this paper.

\bibliographystyle{plain}
\bibliography{papers,internet}

\end{document}